\documentclass[12pt, a4paper, reqno]{amsart}

\usepackage[english]{babel}
\usepackage[T1]{fontenc}
\usepackage[utf8]{inputenc}

\usepackage{lmodern}
\usepackage{amsmath}
\usepackage{amssymb}
\usepackage{amsthm}
\usepackage{url}
\usepackage{xcolor}

\usepackage{hyperref} 

\usepackage{enumerate}
\usepackage[normalem]{ulem}

\usepackage{float}
\restylefloat{table}

\renewcommand{\epsilon}{\varepsilon}
\renewcommand{\theta}[0]{\vartheta}
\renewcommand{\phi}[0]{\varphi}

\newcommand{\Z}{\mathbb{Z}}
\newcommand{\Q}{\mathbb{Q}}

\newcommand{\Size}[1]{\left\lvert #1 \right\rvert}

\newcommand{\Span}[1]{\left\langle\, #1 \,\right\rangle}

\newcommand{\Set}[1]{\left\{ #1 \right\}}

\newcommand{\norm}[0]{\trianglelefteq}

\DeclareMathOperator{\End}{End}
\DeclareMathOperator{\Aut}{Aut}

\DeclareMathOperator{\Hol}{Hol}

\newtheorem{dummy}{Dummy}
\numberwithin{dummy}{section}
\numberwithin{figure}{section}

\newtheorem{theorem}[dummy]{Theorem}

\newtheorem*{thestars}{Theorem}
\newtheorem*{thea}{Theorem A}

\newtheorem{lemma}[dummy]{Lemma}

\newtheorem{prop}[dummy]{Proposition}

\newtheorem{corollary}[dummy]{Corollary} 

\theoremstyle{definition}
\newtheorem{definition}[dummy]{Definition}

\theoremstyle{remark}

\newtheorem{remark}[dummy]{Remark}

\newtheorem{fact}[dummy]{Fact}

\numberwithin{equation}{section}

\newcommand{\Gal}{\mathop{\mathrm{Gal}}\nolimits}

\begin{document}

\title[On the ranks of a brace]{On the ranks of the additive and\\
      the multiplicative groups of a brace}


\author{A. Caranti}

\address[A.~Caranti]%
 {Dipartimento di Matematica\\
  Universit\`a degli Studi di Trento\\
  via Sommarive 14\\
  I-38123 Trento\\
  Italy} 

\email{andrea.caranti@unitn.it} 

\urladdr{http://science.unitn.it/$\sim$caranti/}

\author{I. Del Corso}

\address[I.~Del Corso]%
        {Dipartimento di Matematica\\
          Universit\`a di Pisa\\
          Largo Bruno Pontecorvo, 5\\
          56127 Pisa\\
          Italy}
\email{ilaria.delcorso@unipi.it}

\urladdr{http://people.dm.unipi.it/delcorso/}

\subjclass[2010]{16T25 20D15 20B35 20E18}

\keywords{braces, skew braces, Yang-Baxter equation, holomorph,
          regular subgroups, 
          Hopf-Galois structures}   

\begin{abstract}

In    \cite[Theorem    2.5]{Bac16}    Bachiller   proved    that    if
$(G, \cdot, \circ)$ is  a brace of order the power of  a prime $p$ and
the rank  of $(G,\cdot)$ is  smaller than $p-1$,  then the order  of any
element is  the same  in the additive  and multiplicative  group. This
means that in this case the isomorphism type of $(G,\circ)$ determines
the isomorphism type of $(G,\cdot)$.

In this paper  we complement Bachiller's result in  two directions. In
 Theorem~\ref{teo:small_rank} we prove that  if $(G, \cdot, \circ)$ is
 a  brace  of  order  the  power of  a  prime  $p$,  then  $(G,\cdot)$
 has  \emph{small}  rank   (\emph{i.e.}  $<  p-1$)  if   and  only  if
 $(G,\circ)$  has  \emph{small} rank.   We  also  provide examples  of
 groups of rank $p-1$ in which  elements of arbitrarily large order in
 the  additive  group become  of  prime  order in  the  multiplicative
 group. When the rank is larger, orders may increase.
\end{abstract}

\thanks{The     authors   are   members   of
  INdAM---GNSAGA. The authors gratefully  acknowledge support from the
  Departments of  Mathematics of  the Universities  of 
  Pisa and Trento.  The second  author has performed this activity in
    the framework of  the PRIN 2017, title  ``Geometric, algebraic and
    analytic methods in arithmetic''.}

\maketitle

\thispagestyle{empty}

\section{Introduction}

Let  $L/K$ be  a  finite field  extension,  and let  $H$  be a  finite
cocommutative $K$-Hopf  algebra. $H$ defines a  Hopf-Galois structure
on  $L/K$ if  there  exist a  $K$-linear  map $\mu\colon  H\mapsto{\rm
  End}_K(L)$  giving  $L$  a  left $H$-module  algebra  structure  and
inducing  a  $K$-vector  space  isomorphism  $K\otimes_K  H\mapsto{\rm
  End}_K(L)$.  This notion  was introduced  by Chase  and Sweedler  in
\cite{CS}.  Greither and Pareigis in~\cite{GP} showed that finding the
Hopf-Galois structures can be reduced to a group-theoretic problem.

In the particular case when $L/K$ is a Galois extension we can state 
Greither and Pareigis results as follows.

\begin{thestars}
  Let $L/K$ be a finite Galois extension of fields
  and  let $\Gamma=\Gal(L/K)$.   There is  a bijective  correspondence
  between the set of isomorphism  classes of Hopf-Galois structures on
  $L/K$ and the set of regular  subgroups $G$ of the group $S(\Gamma)$
  of permutations  on the  set $\Gamma$, which  are normalised  by the
  image $\rho(\Gamma)$  of the right regular  representation $\rho$ of
  $\Gamma$. 
\end{thestars}

(In this paper we  use the  right regular  representation $\rho$.  In the
literature it  is more common  to use the left  regular representation
$\lambda$.)

The groups  $G$ and $\Gamma$ have  the same cardinality but  they need
not be  isomorphic.  A Hopf-Galois structure $H$ is called of
\emph{type}  $G$,  if  $G$  is  the  group  associated  to  $H$ in  the
Greither-Pareigis correspondence.

Childs~\cite{Chi89} and Byott~\cite{Byo96} observed that the condition
that $\rho(\Gamma)$ normalises $G$ can  be reformulated by saying that
$\Gamma$ is contained in the holomorph $\Hol(G)$ of $G$, regarded as a
subgroup of the group $S(G)$ of permutations on the set $G$, an
advantage being that $\Hol(G)$ is usually much smaller than
$S(\Gamma)$.  Therefore, the number of 
Hopf-Galois  structures on $L/K$ of  type $G$ can be  computed in
terms of
the number of regular subgroups  of the holomorph $\Hol(G)$ of $G$,
which are
isomorphic  to  $\Gamma$  (see~\cite[Corollary  p.~3320]{Byo96}).  In
particular,  $L/K$ admits  a structure  of type  $G$  precisely when  the
holomorph $\Hol(G)$ of  $G$ contains a regular  subgroup isomorphic to
$\Gamma$.

The study of Hopf-Galois structures,  or equivalently, of regular subgroups of
holomorphs, is strictly related to  the  theory of \emph{skew (left)
  braces}. In  fact, if $G$ is  a group with respect  to the operation
``$\cdot$'',  classifying  the  regular   subgroups  of  $\Hol(G)$  is
equivalent to determining the operations  ``$\circ$'' on $G$ such that
$(G, \cdot, \circ)$ is a (right) skew brace~\cite{skew}, that is, $(G,
\circ)$ is also a  group, and the two group structures  on the set $G$
are linked by the identity
\begin{equation}
  \label{eq:brace}
  (g \cdot h) \circ k = (g \circ k) \cdot k^{-1} \cdot (h \circ k),
\end{equation}
for all $g, h, k \in G$.
This connection was first  observed by Bachiller in~\cite[\S 2]{Bac16}
and it is described in detail in the appendix to~\cite{SV2018}. If
$(G, \cdot, \circ)$ is a skew brace, $(G, \cdot)$ is its
\emph{additive} group, and $(G, \circ)$ its \emph{multiplicative}
group. A 
\emph{(left) brace} can be defined as a skew brace with a commutative
additive group, but the theory of braces predates that of skew
braces~\cite{braces}.

In recent years, these different approaches concurred to construct a
rich theory.  
A number of papers are devoted to enumerating Hopf-Galois structures
on Galois extensions of degree of a particular form
(\cite{Byo96, kohl98, Byott04, Childs05, zen2018, AlBy20-sqrf, AB2020-pq,
  AB2020-p2q, CCDC}). 
Part of the literature is devoted to understanding the group-theoretic
relation between $\Gamma$ and $G$ for a Galois extension with Galois
group isomorphic to $\Gamma$ could admit a Hopf-Galois structure of
type $G$ or not (\cite{FCC}, \cite{Byott13}, \cite{Byott15}, \cite{Bac16}, \cite{Tsang19}, \cite{Nasy19}, \cite{TQ20}). 

In the language of skew braces, this correspond to understanding how
properties of $(G, \cdot)$ influence those of $(G,\circ)$. and vice versa,
when $(G,\cdot,\circ)$ is a skew brace.

In  \cite[Theorem 2.5]{Bac16} Bachiller proved that if 
$(G,\cdot,\circ)$ is  a brace of order the power of a prime $p$ and the rank of $(G,\cdot)$ is lower than $p-1$, then
the  order   of  any  element  is   the  same  in  the   additive  and
multiplicative group. This means that in this case  the
isomorphism type  of $(G, \circ)$  determines the isomorphism  type of
$(G, \cdot)$. 
Bachiller's Theorem 2.5 generalises~\cite[Theorem 1]{FCC} using a similar argument.

In this paper we adopt the same point of view as in \cite{FCC} and \cite{Bac16} and study some relations between the additive and the multiplicative group of a brace, focusing on the order of the elements. Our first result is Theorem~\ref{teo:small_rank} in which we prove that the rank of $(G,\cdot)$ is lower than $p-1$ if and only if  the rank of $(G,\circ)$ is lower than $p-1$ (see Definition~\ref{def:small_rank}). This result builds upon \cite[Theorem2.5]{Bac16}.

When   the  rank   of  $(G,   \cdot)$  is   $p  -   1$,  we   show  in
Proposition~\ref{prop:does-not-increase} that  the orders  of elements
may only  decrease when going  from $(G,  \cdot)$ to $(G,  \circ)$; in
Proposition~\ref{prop:examples} we provide  examples in which elements
of arbitrarily large order in the additive group become of prime order
in the  multiplicative group.  When the  rank is larger than  $p - 1$,
orders of elements may also increase. 

 The paper is enriched with some corollaries in which we specify the consequences of our result in the Hopf-Galois context.

Section~\ref{sec:statements} contains  the statements of  our results,
and       
a   description    of    the   method of the gamma
function (see~\cite{CCDC})   we   use.     In 
Section~\ref{sec:FCC}   we   prove  Theorem~\ref{teo:small_rank}   and
Proposition~\ref{prop:does-not-increase}.                           In
Section~\ref{sec:examples} we prove Proposition~\ref{prop:examples}.

\section{Statements}
\label{sec:statements}

Given a group $G = (G, \cdot)$, write $S(G)$ for the group of
permutations on the set $G$. Let $\rho$ be the right regular
representation of $G$,
\begin{align*}
  \rho :\ &G \to S(G)\\
          &g \mapsto (x \mapsto x \cdot g)
\end{align*}
The (permutational) holomorph $\Hol(G)$ is the normaliser
\begin{equation*}
  \Hol(G)
  =
  N_{S(G)} (\rho(G))
\end{equation*}
of the image $\rho(G)$ of $\rho$ in $S(G)$. $\Hol(G)$ is isomorphic to the
(abstract) holomorph $\Aut(G) \ltimes G$ of $G$.

In the holomorph of an abelian $p$-group $G$ of order at least $p^{3}$
one can find (non-abelian) regular  subgroups which are not isomorphic
to $G$.   In the other direction,  by \cite[Theorem 2.5]{Bac16} (quoted below as Theorem A), 
when the group  $G$ is \emph{small}, then the
isomorphism class  of a regular  subgroup of $\Hol(G)$  does determine
the isomorphism class of $G$.

Let $N$ be a regular subgroup of $S(G)$. The map $\nu : G \to N$, that
takes  $g \in  G$ to  the  unique element  $\nu(g)$ of  $N$ such  that
$1^{\nu(g)} = g$ is a bijection.  (We write the action of permutations
as exponents.)  Using $\nu$ for transport of structure from $N$ to $G$
yields a group operation ``$\circ$'' on  $G$ such that $\nu(g \circ h)
= \nu(g) \nu(h)$, so that 
\begin{equation*}
  \nu : (G, \circ) \to N
\end{equation*}
is an
isomorphism. Moreover $g^{\nu(h)} = g \circ h$, for $g, h \in G$. 

As in~\cite[Section~2]{CCDC},  to which we refer  for further details,
we have that the regular subgroup $N$ normalises $\rho(G)$ if and only
if there is a function $\gamma : G \to \Aut(G)$ such that
\begin{equation}
  \label{eq:GF}
  \nu(g) = \gamma(g) \rho(g),
  \quad
  \text{for $g \in G$}.
\end{equation}
We thus have
\begin{equation*}
  h \circ g
  =
  h^{\nu(g)}
  =
  h^{\gamma(g) \rho(g)}
  =
  h^{\gamma(g)} \cdot g.
\end{equation*}

The functions $\gamma : G \to \Aut(G)$ appearing in~\eqref{eq:GF} are
characterised by the functional equation
\begin{equation*}
  \gamma(h^{\gamma(g)} \cdot g)
  =
  \gamma(h) \gamma(g)
  \quad
  \text{for $h, g \in G$}.
\end{equation*}
As in~\cite{CCDC},  we call  this equation the  \emph{gamma functional
  equation},  and refer  to the  functions
$\gamma$ that satisfy  it as \emph{gamma functions}.

A (right) skew brace \cite{GV} is a triple $(G, \cdot, \circ)$, where $G$ is a
set, ``$\cdot$'' and ``$\circ$'' are two group operations on $G$, and
the following brace axiom holds for $a, b, c \in G$
\begin{equation}
  \label{eq:brace-axiom}
  ((a \cdot b) \circ c)) \cdot c^{-1}
  =
  (a \circ c) \cdot c^{-1}
  \cdot
  (b \circ c) \cdot c^{-1}.
\end{equation}

If $(G, \cdot, \circ)$ is a skew brace, then~\eqref{eq:brace-axiom}
and the fact that $(G, \circ)$ is a group yield that for all $c \in G$ the
maps $G \to G$ given by 
\begin{equation*}
  \gamma(c)
  :
  a \mapsto (a \circ c) \cdot c^{-1}
\end{equation*}
are automorphisms of $G$, and that $\gamma : G \to \Aut(G)$ is a gamma
function. Conversely, if $\gamma : G \to \Aut(G)$ is a gamma
function, then $(G, \cdot, \circ)$ is a skew brace, where $a \circ b =
a^{\gamma(b)} \cdot b$, for $a, b \in G$. Moreover, the set $N = \Set{
  \nu(g) : g \in G}$ of the functions
\begin{align*}
  \nu(g) :\ &G \to G\\
            &x \to x \circ g
\end{align*}
is a regular subgroup of $\Hol(G)$ isomorphic to $(G, \circ)$.
(See for instance the discussion after Theorem 2.2 in~\cite{CCDC}.)

In the following, given a group $G = (G, \cdot)$, we will make use
without further mention of the equivalence described in this section
among the following concepts:
\begin{enumerate}
\item
  a regular subgroup $N$ of $\Hol(G)$,
\item
  a gamma function on $G$,
\item
  a skew brace with additive group $(G, \cdot)$.
\end{enumerate}

Note that braces predate skew  braces (\cite{braces, GV}), but a brace
may be  defined as a skew  brace whose additive group  is abelian.

If
$(G, \cdot, \circ)$ is a (skew) brace, we refer to $(G, \cdot)$ as its
\emph{additive group}, and to $(G, \circ)$ as its \emph{multiplicative
  group}.
  
  In the paper \cite{Bac16} Bachiller, generalizing the main result of \cite{FCC}, proved 
  that if
$(G,\cdot,\circ)$ is  a brace of order the power of a prime  $p$, and the rank of the abelian group $(G,\cdot)$ is less than $p-1$
 then the
isomorphism type  of $(G, \circ)$  determines the isomorphism  type of
$(G, \cdot)$.  This result depends on  the fact that in  these braces,
the  order   of  any  element  is   the  same  in  the   additive  and
multiplicative group. 
  
  \begin{thea}[\protect{\cite[Theorem 2.5]{Bac16}}]
  Let $p$ be a prime, and 
  let $(G, \cdot, \circ)$ be a brace of order a power of the prime
  $p$, with $(G, 
  \cdot)$ of rank $\null < p - 1$.
  
  Then  each element has the same order in $(G, \cdot)$ and $(G, \circ)$.
  
  Moreover, if $(G, \cdot)$ is abelian,then  $(G, \cdot)$ and $(G, \circ)$ are isomorphic.
  
\end{thea}

In her recent paper \cite{Crespo}, Teresa Crespo considers
the same question in the Hopf Galois context.  In
the     language of braces  her   result  states that if
$(G,\cdot,\circ)$ is  a brace of order  $p^{n}$ where $p \ge 3$ is a
prime  and $n<p$, then the
isomorphism type  of $(G, \circ)$  determines the isomorphism  type of
$(G, \cdot)$. Her result  \cite[Theorem  6]{Crespo} improves Theorem A in the case when $p$ is odd and $(G,\cdot)$ is elementary abelian.

In this paper we complement Theorem A in two directions. In Theorem~\ref{teo:small_rank} we give a sort of converse of Theorem A, and in Proposition~\ref{prop:does-not-increase} we consider the case when the rank of $(G,\cdot)$ is $p-1$, improving \cite[Theorem  6]{Crespo}.

\begin{definition}
  \label{def:small_rank}
  Let $\mathcal{G}$ be a finite $p$-group, for a prime $p$. Its
  \emph{rank} $r_{p}$  is the
  maximum $r$ such that $\mathcal{G}$  has a subgroup of exponent $p$ and
  order $p^{r}$.
  We say  that $\mathcal G$ has  \emph{small rank} if $r_{p} < p - 1$ 
   
  Let $\mathcal{G}$ be a finite group, and let $p$ be a prime dividing
  its order. We say that has  \emph{small p-rank} if  a Sylow
  $p$-subgroup of $\mathcal{G}$ has  small rank.
  We say  that $\mathcal{G}$ has
  \emph{small  rank} if  it has  a small  $p$-rank for  each prime  $p$
  dividing its order.
\end{definition}

The following theorem will be proved in
Section~\ref{sec:FCC}. 
\begin{theorem}
  \label{teo:small_rank}
  Let $p$ be a prime, and 
  let $(G, \cdot, \circ)$ be a brace of order a power of the prime
  $p$.

  Then 
  $(G,\cdot)$ has small rank if and only if 
  $(G,\circ)$ has small rank.

  When these conditions hold,
  $(G, \cdot)$ and $(G, \circ)$
  have   the same rank, and the same number of elements of each order.

\end{theorem}

\begin{remark}
  \label{rem:Berkovich}
  Yakov Berkovich has  shown in~\cite[Proposition~1.6(a)]{Ber00} (see
  also  the   paragraph  preceding   Proposition~7.8  of~\cite{Ber02},
  and~\cite[Lemma 3(a)]{Ber05}) that a  finite $p$-group of small rank
  is regular in the sense of Philip Hall~\cite{Hall}.
\end{remark}
The  following  corollaries  of Theorem~\ref{teo:small_rank}  will  be
proved in Section~\ref{sec:proofs_cor}.

The first one generalises  \cite[Theorem 6]{Crespo} for the case of groups of rank  $<p-1$ and is indeed a consequence of \cite[Thoerem 2.5]{Bac16}.

\begin{corollary}
  \label{cor:A_plus}
  Let $p$ be a prime, and  let $G_{1}, G_{2}$ be two abelian groups of
  order the same power of a prime $p$, with $G_{1}$ of small rank.

  Let $N_{i}$ be a regular subgroup of the
  holomorph $\Hol(G_{i})$ of $G_{i}$, for $i = 1, 2$.

  If $N_{1} \cong N_{2}$, then $G_{1} \cong G_{2}$. In particular,
  $G_{2}$ has also small rank.
\end{corollary}

\begin{corollary}
  \label{cor:pochi_pelementi}
  Let $L/K$ be a Galois extension of degree  a power of a prime $p$, let
  $\Gal (L/K) = \Gamma$, and let $r$ denotes the rank of $\Gamma$. Then
  \begin{enumerate}
  \item
    \label{item:corppe1}
    If $r < {p-1}$ (that is, $\Gamma$ has small rank)
    and  $L/K$ admits abelian Hopf-Galois structures, then all of them are of the same type and the
    group $G$ associated to these structures  is determined by $\Gamma$
    and has the same rank as $\Gamma$. 
    In particular, if $\Gamma$ is abelian, then  every abelian Hopf Galois structure on $L/K$ is of type $\Gamma$.
      \item
    \label{item:corppe2}
    If  $r \ge {p-1}$  (that is,  $\Gamma$ does  \emph{not} have  small
    rank), then every abelian group  giving a Hopf-Galois structure on
    $L/K$ has rank $\null\ge{p-1}$.
  \end{enumerate}
\end{corollary}

\begin{remark}
\label{rem:not_any_p_group}
  Corollary~ \ref{cor:pochi_pelementi} delimits the types of the possible abelian Hopf-Galois structure on a Galois extension of prime power order, dividing the abelian $p$-groups in two rigid classes, accordingly with their rank.
  
  However, it gives no indication on the question of whether
  a Galois extension with a non-abelian Galois group $\Gamma$ does admit at least an abelian Hopf-Galois structure. This is the same as asking whether the (non-abelian) groups $\Gamma$ can be
  the  multiplicative group  of  a brace.   A  necessary condition  is
  provided by Byott, who showed in~\cite[Theorem~1]{Byott15}
  that for $L/K$ could admit a nilpotent Hopf-Galois structure its Galois Group $\Gamma$ must be solvable.
 In particular this implies that  the
  multiplicative group of a brace is solvable (see also \cite[Theorem 2.15]{EST99}). 
  
  In the paper \cite{Bac16} Bachiller considered the converse, asking whether any finite solvable group is the multiplicative group of a brace, giving to this question a negative answer. In fact, he provided an example of a $p$-group $\Gamma_0$ of order $p^{10}$ with all elements of order $p$ which, for some large value of the prime $p$, is not the multiplicative group of a brace.
Therefore, the Galois extensions with Galois group isomorphic to  $\Gamma_0$ do not admit abelian Hopf-Galois structures. Since all $p$-groups are realisable as Galois group over $\Q$, then this actually happens also in the number field context. 
 \end{remark}

\begin{remark}
  In \cite[Remark 9]{Crespo} Crespo points out that a Galois extension
  with  Galois  group  $C_9\times   C_3\times  C_3$  has  Hopf  Galois
  structure of types $C_9^2$ and $C_3^4$, in addition to the classical
  one.    This     shows    that     under    the     hypothesis    of
  Corollary~\ref{cor:pochi_pelementi}.\eqref{item:corppe2},  a  Galois
  extension can have  abelian Hopf-Galois structures of different types.
\end{remark}

\begin{corollary}
  \label{cor:orders-gen}
  Let  $(G, \cdot, \circ)$ be a brace of finite order.

  Let $p$ be a prime divisor of the order $G$.
  
  Then the
  $p$-rank of $(G,\cdot)$ is small (that is, $\null<p-1$) if and only if
  the $p$-rank of $(G, \circ)$ is small. In this case
  the Sylow $p$-subgroups of $(G, \circ)$ and $(G, \cdot)$ have the same
  number of elements of each order.
  
  In particular, $(G, \cdot)$ has small rank if and only if $(G, \circ)$
  has  small rank. If  this  is the  case, and  $(G, \circ)$  is
  abelian then $(G, \cdot)\cong(G, \circ)$.
\end{corollary}

\begin{corollary}
  \label{cor:esp_piccoli}
  Let  $L/K$  be a finite Galois  extension  and  let $\Gal  (L/K) =
  \Gamma$.  
  
  If $\Gamma$  has  small  rank,  
  and  $L/K$ admits abelian Hopf-Galois structures, then all of them are of the same type $G$, and  the group  $G$   is determined by $\Gamma$. 
    
      In particular, if $\Gamma$ is abelian, then every abelian Hopf Galois structure on $L/K$ is of type $\Gamma$.
\end{corollary}

The   following   proposition   shows   that  the   second   part   of
Theorem~\ref{teo:small_rank} fails  when $(G,\cdot)$  (or equivalently
$(G,\circ)$) does not have small rank.

\begin{prop}
  \label{prop:examples}
  For
  every $k > 1$, there is a brace $(G, \cdot, \circ)$ of $p$-power order,
  with $(G, \cdot)$ of
  rank $p-1$, with the following properties.
  \begin{enumerate}
  \item
    $(G, \circ)$ is non-abelian.
  \item
    There is a maximal subgroup $H$ of $(G, \cdot)$, such that every element of
    $H$ has the same order in $(G, \cdot)$ and $(G, \circ)$.
  \item
    Every  element $g \in G \setminus H$  has order $p^{k}$ in $(G,
    \cdot)$, and order
    $p$ in $(G, \circ)$.
  \end{enumerate}
\end{prop}

In the  examples of  Proposition~\ref{prop:examples}, the order  of an
element  does  not increase  when  going  from  $(G, \cdot)$  to  $(G,
\circ)$.  This is  actually a general fact, as shown  in the following
\begin{prop}
  \label{prop:does-not-increase}
  Let $(G, \cdot, \circ)$ be a brace  of $p$-power order,
  with $(G, \cdot)$ of
  rank $p-1$.

  Then the order of an
  element does not increase when going from $(G, \cdot)$ to $(G,
  \circ)$.
\end{prop}
In the case when $(G, \cdot)$ is elementary abelian the previous proposition gives the following corollary, which was already covered in \cite{Crespo}. 

\begin{corollary}
\label{cor:teob}
Let $(G, \cdot, \circ)$ be a brace  of $p$-power order,
  with $(G, \cdot)$ isomorphic to $C_p^{m}$ with $m\le p-1$. Then each element has the same order in $(G,\cdot)$ and  $(G, \circ)$.
In particular  if $(G,\circ)$ is abelian, then it is isomorphic to $(G, \cdot)$.
\end{corollary}

When the rank of $G$ reaches $p$, the order of an element may increase
when going  from $(G, \cdot)$ to  $(G, \circ)$, as shown  for instance
by~\cite[Example~8]{FCC} .
 
The braces $(G, \cdot,  \circ)$ of Proposition~\ref{prop:examples} are
bi-skew  braces (\cite{Childs-bi,  Caranti-bi}), so  that $(G,  \circ,
\cdot)$ is a skew brace (with non-abelian additive group $(G, \circ)$)
in which  the order of  an element may  increase when going  from $(G,
\circ)$ to $(G, \cdot)$.

\section{Proofs of Theorem~\ref{teo:small_rank} and
Proposition~\ref{prop:does-not-increase}}

\label{sec:FCC}

For a  finite $p$-group  $(H, *)$, and  $i \ge 0$,  we will  denote by
$\Omega_{i}(H, *)$  the \emph{set}  of elements of  $(H, *)$  of order
dividing $p^{i}$.  

We begin with proving Proposition~\ref{prop:does-not-increase}, which  is the particular case $m=p-1$ of 
the following lemma.
\begin{lemma}
  \label{lemma:lemma}
    
  Let $p$ be a prime.
  
  Let $G = (G, \cdot)$ be a finite abelian $p$-group of $p$-rank $m$.

  Let $(G, \cdot, \circ)$ be a brace.

    If $m \le p - 1$, then for each $i \ge 0$ one has
    \begin{equation}
      \label{eq:a_tale_of_one_Omega}
      \Omega_{i+1}(G, \cdot)
      \subseteq
      \Omega_{i+1}(G, \circ).
    \end{equation}
    In other words, the order of an element does not increase when going
    from $(G, \cdot)$ to $(G, \circ)$.
\end{lemma}

Note that in  the Lemma the $\Omega_{i}(G,  \cdot)$ are characteristic
subgroups of  the abelian group $(G,  \cdot)$. On the other  hand, the
$\Omega_{i}(G, \circ)$ need  not be subgroups of $(G,  \circ)$; in the
examples of Section~\ref{sec:examples}, no $\Omega_{i}(G, \circ)$ is a
subgroup of $(G,  \circ)$, for $0 < i <  \log_{p}(\exp(G))$.  In these
examples  we  have  $m  =  p  -   1$;  however,  when  $m  <  p  -  1$
Remark~\ref{rem:Berkovich} implies that the $\Omega_{i}(G, \circ)$ are
also subgroups of $(G, \circ)$.

\begin{proof}[Proof of Lemma~\ref{lemma:lemma}]
  For $g \in G$, write
  \begin{equation*}
    \delta(g) = - 1 + \gamma(g) \in \End(G, \cdot).
  \end{equation*}

  It is immediate to see that for $g \in G$ we have, for the $p$-th
  power $g^{\circ p}$ of $g$ in $(G, \circ)$,
  \begin{equation*}
    g^{\circ p}
    =
    g^{\gamma(g)^{p-1} + \dots + \gamma(g) + 1},
  \end{equation*}
  from which we obtain the formula
  \begin{equation}
    \label{eq:powp}
    g^{\circ p}
    =
    g^{p
      +
      \binom{p}{2} \delta(g)
      + \dots +
      \binom{p}{p-1} \delta(g)^{p-2}}
    \cdot
    g^{\delta(g)^{p-1}}.
  \end{equation}
  (We write simply $g^{p}, g^{\binom{p}{2}}$, etc.\ for the powers of $g$ in
  $(G, \cdot)$.) 

  We first assume $m \le p - 1$, and prove~\eqref{eq:a_tale_of_one_Omega}.
  
  Let   us start  with the  case $i  = 0$.   Let $g  \in \Omega_{1}(G,
  \cdot)$.    Then  in~\eqref{eq:powp}   we   have   $g^{\circ  p}   =
  g^{\delta(g)^{p-1}}$. Now $\gamma(G) \Omega_{1}(G,
    \cdot)$   is  a   finite   $p$-group,   thereby  nilpotent,   with
    $\Omega_{1}(G, \cdot)$ a  normal subgroup of it.  Therefore for all $i \ge 1$ we have
    \begin{itemize}
    \item 
      either $[\Omega_{1}(G,  \cdot),
        \underbrace{\gamma(G),  \dots,  \gamma(G)}_{i-1}] = 1$,
    \item
      or $[\Omega_{1}(G, \cdot),
        \underbrace{\gamma(G),     \dots,    \gamma(G)}_{i-1}]     >
      [\Omega_{1}(G,    \cdot),     \underbrace{\gamma(G),    \dots,
          \gamma(G)}_{i}]$.
    \end{itemize}
   Since $\Omega_{1}(G,    \cdot)$ has  order at most
    $p^{p-1}$, we obtain
    \begin{equation*}
      [\Omega_{1}(G,     \cdot),     \underbrace{\gamma(G),     \dots,
          \gamma(G)}_{p-1}] = 1.
  \end{equation*}
   
  Now note that for $h, g \in G$. one has that $h^{\delta(g)}$ equals the
  commutator $[h, \gamma(g)]$ in the (abstract) holomorph of $G$.
  We thus obtain $g^{\delta(g)^{p-1}} = 1$, and we are done.

  Proceeding by induction, we take $i \ge 1$, and assume
  \begin{equation*}
    \Omega_{i}(G, \cdot)
    \subseteq
    \Omega_{i}(G, \circ),
  \end{equation*}
  and we prove
  \begin{equation*}
    \Omega_{i+1}(G, \cdot)
    \subseteq
    \Omega_{i+1}(G, \circ).
  \end{equation*}

  Let   $g   \in   \Omega_{i+1}(G,  \cdot)$. We have
  $g^{p}  \in \Omega_{i}(G,  \cdot)$, so  that by~\eqref{eq:powp}  the
  following are equivalent
  \begin{enumerate}
  \item
    $g^{\circ p} \in \Omega_{i}(G, \cdot)$, and
  \item
    $g^{\delta(g)^{p-1}} \in \Omega_{i}(G, \cdot)$.
  \end{enumerate}
  Now  $M  = \Omega_{i+1}(G,  \cdot)  /  \Omega_{i}(G, \cdot)$  is  an
  elementary abelian section of $G$, invariant under automorphisms, of
  order  at  most  $p^{p-1}$, so  $g^{\delta(g)^{p-1}}  \in  \Omega_{i}(G,
  \cdot)$. We
  have  obtained  that $g^{\circ  p}  \in  \Omega_{i}(G,  \cdot)  \subseteq
  \Omega_{i}(G, \circ)$, so that $g \in \Omega_{i+1}(G, \circ)$.
\end{proof}

We recall the following observation from~\cite{CCDC}

\begin{lemma}[\protect{\cite[Proposition 2.6]{CCDC}}]
  \label{lemma:2.6}
  Let $G=(G, \cdot)$ be a finite group, let $H \subseteq G$ and let
  $\gamma$ be a 
  GF on $G$.
  
  Any  two of  the
  following conditions imply the third one:
  \begin{enumerate}
  \item\label{i1} $H \le G$;
  \item\label{i2} $(H,\circ) \le (G, \circ)$;
  \item\label{i3} $H$ is $\gamma(H)$-invariant.
  \end{enumerate}
  If these conditions hold, then $(H,\circ)$ is isomorphic to a regular
  subgroup of $\Hol(H).$ 
\end{lemma}

We turn now  to the proof of  Theorem~\ref{teo:small_rank}.
As we noted above, the ``if" part of  this theorem  is \cite[Theorem 2.5]{Bac16}.

Conversely, we have to prove that if $(G,  \circ)$  has  small  rank, then

$(G,  \cdot)$ has also small rank.
 
Consider   the  finite   $p$-group  $K   =  \gamma(G)
  G$. Since $\Omega_{1}(G)$  is a characteristic subgroup  of $G$, and
  $G$ is a  normal subgroup of $K$, it follows  that $\Omega_{1}(G)$ is a
  normal subgroup of $K$.  Since  $K$ is nilpotent, and $\Omega_{1}(G)
  \norm  K$, there  is  a central  series of  $K$  which goes  through
  $\Omega_{1}(G)$.  Refining this  series  to a  principal series,  we
  obtain that for  each divisor of the order  of $\Omega_{1}(G)$ there
  is  a subgroup of $\Omega_{1}(G)$  which is normal in  $K$, and
  thus in particular invariant under $\gamma(G)$.  Thus if the rank of
  $(G, \cdot)$ is greater than or equal to $p-1$, so that the order of
  $\Omega_{1}(G)$ is  greater than or  equal to $p^{p-1}$,  then there
  will be  a subgroup  $T$ of $\Omega_{1}(G)$  of order  $p^{p-1}$ and
  exponent $p$, invariant under $\gamma(G)$.

Now,  Lemma~\ref{lemma:2.6}  ensures that $(T, \circ)$ is a subgroup of $(G, \circ)$,
  and by  Proposition~\ref{prop:does-not-increase}  of Section~\ref{sec:statements},  $(T, \circ)$
  has exponent $p$. This gives a contradiction.

\section{Proofs of the corollaries}
\label{sec:proofs_cor}

\begin{proof}[Proof of Corollary~\ref{cor:A_plus}]
  A regular subgroup $N_i$ of  $\Hol(G_i)$ corresponds to a brace $(G_i,
  \cdot,  \circ)$ with  additive group  $G_i$ and  multiplicative group
  $(G_i,\circ)\cong N_i$.  By Theorem~\ref{teo:small_rank}, if  the rank
  of $G_1$ is $\null< p-1$ then the group $N_1$ determines the structure
  of $G_1$. Therefore, $N_{1} \cong  N_{2}$ implies $G_{1} \cong G_{2}$,
  and clearly the rank of $G_2$ is the same as the rank of $G_1$.
\end{proof}

\begin{proof}[Proof of Corollary~\ref{cor:pochi_pelementi}]
Let $G$ be an abelian group defining a Hopf Galois structure on $L/K$. Then we can 
efine on $G$   an operation  $\circ$  such  that $\Gamma\cong(G,\circ)$  and
  $(G,\cdot,\circ)$   is   a   brace.  The   corollary   follows   from
  Theorem~\ref{teo:small_rank} applied to $(G,\cdot,\circ)$.
\end{proof}

For   the  next   proof,  note   that  Lemma~\ref{lemma:2.6}   can  be
reformulated as follows in the language of skew braces
\begin{corollary}
  \label{cor:lemma2.6}
  Let $ (G, \cdot,  \circ)$ be a skew brace and  let $(H,\cdot)$ be a
  subgroup $(G, \cdot)$.
  
  Then $(H, \cdot, \circ)$ is a  sub-skew brace of $(G, \cdot, \circ)$
  if and only if $(H,\cdot)$ is invariant under $\gamma(H)$.
\end{corollary}

\begin{proof}[Proof of Corollary~\ref{cor:orders-gen}]
  For $p$  a prime divisor of  $\Size{G}$, let $G_p$ denote  the Sylow
  $p$-subgroup  of  $(G,\cdot)$.  Since $G_p$  is  characteristic,  by
  Lemma~\ref{lemma:2.6}        and       Corollary~\ref{cor:lemma2.6},
  $(G_p,\circ)$  is  a subgroup,  actually  a  Sylow $p$-subgroup,  of
  $(G,\circ)$  and  $(G_p,\cdot, \circ)$  is a  sub-brace of  $(G, \cdot,  \circ)$.
  Theorem~\ref{teo:small_rank} applied  to $(G_p,\cdot,  \circ)$ gives
  the  first part,  and applied  to all  prime divisors  of $\Size{G}$
  gives the second one.
\end{proof}

The proofs of Corollaries~\ref{cor:esp_piccoli} and  \ref{cor:teob} are immediate.

\section{Proof of Proposition~\ref{prop:examples}}

\label{sec:examples}

\newcommand{\thegroup}[0]{E}

In this section we construct the examples of
Proposition~\ref{prop:examples}. These are based on the unique pro-$p$
group of maximal class, whose construction we now recall.

Let $p$ be a prime, and $\Z_{p}$ be the ring of $p$-adic integers.
Let $\omega$ be a 
primitive $p$-th root of unity. $\omega$ has minimal polynomial
\begin{equation*}
  x^{p-1} + x^{p-2} + \dots + x + 1 \in \Z_{p}[x]
\end{equation*}
over $\Z_{p}$, so  that the ring $\Z_{p}[\omega]$, when  regarded as a
$\Z_{p}$-module, is free of rank $p - 1$.

The ring $\Z_{p}[\omega]$  is a discrete valuation  ring, with maximal
ideal $I = ( \omega - 1  )$. Consider the automorphism $\alpha$ of the
group  $\thegroup  = (\Z_{p}[\omega],  +)$  given  by multiplication  by
$\omega$. Clearly $\alpha$ has order $p$ in $\Aut(\thegroup)$.

The infinite pro-$p$-group of maximal class is
\begin{equation*}
  M
  =
  \Span{\alpha} \ltimes \thegroup.
\end{equation*}
For $p = 2$ this is the infinite pro-$2$-dihedral group, in which all
elements outside $\thegroup$ have order $2$. In general, we have the following
\begin{fact}
  \label{fact}
  All elements $m \in M \setminus \thegroup$  have order $p$.
\end{fact}
This fact is a statement about the abstract holomorph $\Aut(\thegroup)
\ltimes \thegroup$, and then clearly an analogous fact  holds true  in the
permutational holomorph. 
\begin{proof}
  If $g \in \thegroup$ and $0 < i < p$,  we
  have
  \begin{align*}
    (\alpha^{i} g)^{p}
    &=
    \alpha^{i p} g^{\alpha^{i(p-1)} + \alpha^{i(p-2)} + \dots + \alpha^{i} + 1}
    \\&=
 g^{\alpha^{i(p-1)} + \alpha^{i(p-2)} + \dots + \alpha^i + 1}
    \\&=
    (\omega^{i(p-1)} + \omega^{i(p-2)} + \dots + \omega^i + 1) g
    \\&=
    0
  \end{align*}
  since for $0 < i < p$, $\omega^i$ is a conjugate of $\omega$, that is it has the same minimal polynomial.
    \end{proof}

Consider the group morphism $\gamma : \thegroup \to \Aut(\thegroup)$
which has kernel $I$, and then takes the value $\gamma(1) = \alpha$ on $1$.

Since for $g \in \thegroup$ one has 
\begin{equation*}
  [g, \alpha]
  =
  g^{-1 + \alpha}
  =
  (-1 + \omega) g
  \in I,
\end{equation*}
we have $[\thegroup, \gamma(\thegroup)] = \ker(\gamma)$, and thus
$\gamma$ is a gamma 
function, according to~\cite[Lemma 2.13]{CCDC}.

Such a gamma function $\gamma$ thus defines 
\begin{enumerate}
\item
  a group operation
  \begin{equation*}
    g \circ h = g^{\gamma(h)} \cdot h
  \end{equation*}
  on $\thegroup$  such that $(\thegroup,  \cdot, \circ)$ is a  brace, and
\item
  equivalently a regular subgroup $N$ of $\Hol(\thegroup)$ given by
  \begin{equation*}
    N = \Set{ \gamma(g) \rho(g) : g \in \thegroup }.
  \end{equation*}
\end{enumerate}
Writing $\nu(g) = \gamma(g) \rho(g)$, we have $h^{\nu(g)} = h \circ
g$, and the map
$\nu : (\thegroup, \circ) \to N$ is an isomorphism of groups.

Clearly ``$\cdot$ `` and ``$\circ$'' coincide on $I = \ker(\gamma)$. We
now prove that every element $g \in \thegroup \setminus I$ has order $p$ in
$(\thegroup, \circ)$. Since $\nu$ is an isomorphism of groups, this is
equivalent to showing that all elements $\gamma(g) 
\rho(g)$ of $N$ with $\gamma(g) \ne 1$ have order $p$, and this is
Fact~\ref{fact} above.

The structure of $(\thegroup, \circ)$ is easily seen.
Write $u = 1$ in $\thegroup$ for clarity.
$(I, \circ) \cong (I, +)$ is
an abelian normal 
subgroup of $(\thegroup, \circ)$,
and then for $h \in I$ one has, keeping in mind that
$\gamma(h) = 1$ and $\gamma(u) = \alpha$,
\begin{equation*}
  u^{\ominus 1} \circ h \circ u
  =
  -u^{\gamma(u)^{-1} \gamma(h) \gamma(u)} + h^{\gamma(u)} + u
  =
  h^{\gamma(u)}
  =
  h^{\alpha}
  =
  \omega h.
\end{equation*}
Here $u^{\ominus 1}$ is the inverse of $u$ in $(\thegroup, \circ)$, and we are
using~\cite[Lemma 2.10]{CCDC}.

To finish the proof  of Proposition~\ref{prop:examples}, consider, for
a given $k > 1$, the quotient group
\begin{equation*}
  (G, +) = E / I^{k (p - 1)}.
\end{equation*}
$\alpha$ induces an automorphism of $G$, which we still call
$\alpha$. Write $H = I / I^{k (p - 1)}$. As above, the function $\gamma : G \to
\Aut(G)$ that has kernel $H$, and such that $\gamma(u + I^{k (p - 1)})
= \alpha$, is a gamma function, which defines an operation ``$\circ$''
on $G$. As above, the elements of $G$ have the same order in $(H, +)$
and $(H, \circ)$, while the elements of $G \setminus H$ have order
$p^{k}$ in $(G, +)$, and order $p$ in $(G, \circ)$.

\bibliographystyle{amsalpha}
 
\bibliography{Refs}

\end{document}